\def\be{\begin{eqnarray}}
\def\ee{\end{eqnarray}}

\def\ll#1{\left#1}
\def\r#1{\right#1}
\def\fr{\frac{1}{2}}
\def\mref#1{(\ref{#1})}

\def\bd{\begin{displaymath}}
\def\ed{\end{displaymath}}
\def\k{\mbox{$\kappa$}}
\def\mb{\mbox{}}
\def\ba#1{\begin{array}{#1}}
\def\ea{\end{array}}
\def\nn{\nonumber}
\newfont{\Bbb}{msbm10 scaled 1200}
\documentstyle[12pt,fleqn]{article}

\title{
The contraction of $SU_\mu(2)$\ and its differential structures to $E_\kappa(2)$
}
\author{
Piotr Kosi\'nski$^*$, Pawe\l\ Ma\'slanka\thanks{Supported by KBN grant 2 P03B 130 12}\\
                              Theoretical Physics Department II\\
                                     University of Lodz\\
                          ul. Pomorska 149/153, 90 236 Lodz, Poland\\}
\date{}
\begin{document}
\maketitle
\begin{abstract}
The deformed double covering of $E(2)$\ group, denoted by $\tilde{E}_\kappa(2)$,
is obtained by  contraction  from the $SU_\mu(2)$. The contraction
procedure is then used for producing a new examples of differential calculi:
3D--left covariant calculus on both $\tilde{E}_\kappa(2)$\ and the deformed Euclidean
 group $E_\kappa(2)$\ and two
different 4D--bi\-co\-variant calculi on $\tilde{E}_\kappa(2)$\ which correspond
to the one 4D--bicovariant calculi on $E_\kappa(2)$\ described in Ref.\cite{b14}.
\end{abstract}
\section{Introduction \label{s1}}
The idea to use the contraction procedure in the theory of quantum groups goes
back to E.~Celeghini, R.~Giachetti, E.~Sorace and M.~Tarlini. In a~series of
papers \cite{b1}--\cite{b4} they applied this procedure to quantum deformations
of simple Lie groups producing new examples of quantum groups. Later the contraction
was used by J.~Lukierski, H.~Ruegg, A.~Nowicki and V.~Tolstoy as a tool for obtaining
the \k--deformation of Poincare algebra (\cite{b5},\cite{b6}) and by P.~Zaugg for
constructing the \k--deformation of Poincare group (\cite{b7},\cite{b8}). In \cite{b9}
J.~Sobczyk obtained in this way from $SU_\mu(2)$\ the deformed twodimensional
Euclidean group $E_\kappa(2)$.
In the present paper we use  contraction
as a method for producing a~new examples of noncommutative differential calculi.
For this reason we need a~slightly different contraction technique that the one used in Ref.\cite{b9}.
When applied to $SU_\mu(2)$\ it gives
the deformed counterpart of double covering of $E(2)$, called below $\tilde{E}_\kappa(2)$;\
the projection from $\tilde{E}_\kappa(2)$\ onto $E_\kappa(2)$ can be then easily constructed.
All that is described in sec.\ref{s2}. In sec.\ref{s3} we contract the
3D left--covariant differential calculus on $SU_\mu(2)$\ constructed by S.L.Woronowicz~\cite{b10}.
As a~result we get the 3D left covariant differential calculus on $E_\kappa(2)$.
We also construct the corresponding Lie algebra and prove its equivalence to
$e_\kappa(2)$, the latter being described in Refs.\cite{b11},\cite{b12} as a~dual object
to $E_\kappa(2)$. In sec.\ref{s4} we contract the 4D$_+$\ bicovariant differential calculus
on $SU_\mu(2)$\ described first explicitly by  P.~Stachura,\cite{b13}, and we obtain 4D$_+$\
bicovariant differential calculus on $\tilde{E}_\kappa(2)$\ related to some right
ideal of ``functions" vanishing at the identity of $\tilde{E}_\kappa(2)$.
We give also the corresponding Lie algebra. In sec.\ref{s5} we present the
4D$_-$\ bicovariant differential calculus on $\tilde{E}_\kappa(2)$\ resulting from
the contraction of 4D$_-$\ differential calculus on $SU_\mu(2)$\ (described also in Ref.\cite{b13}).
This  calculus is related to some right
ideal which can be roughly describe as consisting of  ``functions" simultaneously vanishing 
at the identity and some other ``point" 
of $\tilde{E}_\kappa(2)$.
In sec.\ref{s6} we prove that both 4D$_+$\ and 4D$_-$\ calculi on $\tilde{E}_\kappa(2)$\
project onto the same 4D bicovariant differential calculus on $E_\kappa(2)$, first
constructed in Ref.\cite{b14}.

More technical proofs of theorems are relegated to the Appendix.

In the present paper (as well as in all above mentioned papers concerning the contraction procedure) 
the quantum groups are considered in the purely algebraic setting of
deformed universal enveloping algebras.

\section{\k--contraction from $SU_\mu(2)$\ to $E_\kappa(2)$\label{s2}}
In this  section we present a~slightly different scheme of contraction that the
one describe in Ref.\cite{b9}. The contraction procedure presented here allow us to
describe the complete structure of $E_\kappa(2)$\ by considering only the terms
of order $\frac{1}{R}$, where $R$\ is a~contraction parameter.

Let us recall that  $SU_\mu(2)$\  is a~matrix quantum group,\cite{b10}:
\be
&g=\ll(\ba{cc}
\sigma&-\mu\rho^*\\
\rho&\sigma^*\ea\r)\nn&\\
&\Delta(g)=g\hbox{\d{$\otimes$}} g&\nn
\ee
which matrix elements satisfy the following relations:
\be\label{w1}
\rho\rho^*=\rho^*\rho&&\mu(\rho-\rho^*)\sigma=\sigma(\rho-\rho^*)\nn\\
\sigma\sigma^*+\mu^2\rho\rho^*=I&&\mu(\rho+\rho^*)\sigma=\sigma(\rho+\rho^*)\\
\sigma^*\sigma+\rho\rho^*=I&&\nn
\ee
The first step in our scheme is to change  the basis of generators in $SU_\mu(2)$:
\be\label{w2}
\ll(\ba{cc}
\sigma&\mu\rho^*\\
\rho&\sigma^*\ea\r)\to
\fr\ll(\ba{cc}
1&1\\
-1&1
\ea\r)\ll(\ba{cc}
\sigma&-\mu\rho^*\\
\rho&\sigma^*\ea\r)\ll(\ba{cc}
1&-1\\
1&1\ea\r)=\\
=\fr\ll(\ba{cc}
\sigma+\rho-\mu\rho^*+\sigma^*&-\sigma-\rho-\mu\rho^*+\sigma^*\\
-\sigma+\rho+\mu\rho^*+\sigma^*&\sigma-\rho+\mu\rho^*+\sigma^*\ea\r)\nn
\ee

The contraction  ($R$)\ and the deformation (\k) parameters are introduced
by the relations:
\be\label{w3}
\sigma+\rho-\mu\rho^*+\sigma^*&=&a(R)\nn\\
-\sigma-\rho-\mu\rho^*+\sigma^*&=&\frac{w(R)}{R}\\
\mu=e^{\frac{1}{\kappa R}}\nn
\ee

We assume that $\alpha(R)$\ and $w(R)$\ have well--defined limits as $R\to\infty$:
\be\label{w4}
\lim\limits_{R\to\infty}a(R)=a_0&&\lim\limits_{R\to\infty}w(R)=w_0
\ee
The relations \mref{w3} allow us to express the generators  $\rho$\ and $\sigma$\ 
as the functions of the contraction parameter:
\be\label{w5}
\sigma&=&\frac{1}{1+\mu}\ll(\mu a+a^*+\frac{1}{R}(w^*-\mu w)\r)\nn\\
\sigma^*&=&\frac{1}{1+\mu}\ll(\mu a^*+a+\frac{1}{R}(w-\mu w^*)\r)\\
\rho&=&\frac{1}{1+\mu}\ll(a-a^*-\frac{1}{R}(w+w^*)\r)\nn\\
\rho^*&=&-\frac{1}{1+\mu}\ll(a-a^*+\frac{1}{R}(w+w^*)\r)\nn
\ee

We are now ready to perform the contraction of $SU_\mu(2)$\ by taking the limit
$R\to\infty$\ in eq.\mref{w1}. As a~result we obtain the quantum counterpart
of the double covering of $E(2)$, denoted below by $\tilde{E}_\kappa(2)$. The structure
of $\tilde{E}_\kappa(2)$\ is described by the following:
\newtheorem{moje}{Theorem}
\begin{moje}
\mb\newline
$\tilde{E}_\kappa(2)$\ is Hopf algebra generated by the elements: $a_0,\;a_0^*,\;
w_0$\ and $w_0^*$\ subject to the relations:
\be\label{w6}
\mb[a_0^*,w_0]&=&-[a_0^*,w_0^*]=\frac{1}{2\kappa}({a_0^*}^2-I)\nn\\
\mb[a_0,w_0]&=&-[a_0,w_0^*]=\frac{1}{2\kappa}(a_0^2-I)\nn\\
\mb[w_0,w_0^*]&=&-\frac{1}{2\kappa}(a_0+a_0^*)(w_0+w_0^*)\nn\\
\mb[a_0,a_0^*]&=&0\;\;\;\;\;a_0^*a_0=a_0a_0^*=I\nn\\
\Delta a_0&=&a_0\otimes a_0\;\;\;\;\;\Delta a_0^*=a_0^*\otimes a_0^*\\
\Delta w_0&=&w_0\otimes a_0^*+a_0\otimes w_0\nn\\
\Delta w_0^*&=&w_0^*\otimes a_0+a_0^*\otimes w_0^*\nn\\
S(a_0)&=&a_0^*\;\;\;\;\;\;S(a_0^*)=a_0\nn\\
S(w_0)&=&-a_0^*w_0a_0\;\;\;\;\;S(W_0^*)=-a_0w_0^*a_0^*\nn\\
\epsilon(a_0)&=&\epsilon(a_0^*)=1\;\;\;\;\;\epsilon(w_0)=\epsilon(w_0)=0\nn
\ee
\end{moje}

For proof
see Appendix A.

Let us stress that the above structure was obtained by considering only the
terms of order $\frac{1}{R}$\ (cf. the proof of Theorem~1). In order to get the
deformed twodimensional Euclidean group $E_\kappa(2)$\ we put:
\be\label{w7}
A=a_0^2,\;\;\;A^*={a_0^*}^2,\;\;\;v_+=-ia_ow_0,\;\;\;v_-=iw_0^*a_0
\ee
As a~result the following Hopf algebra $E_\kappa(2)$\ (cf. \cite{b11}, \cite{b12},
\cite{b14}) is obtained from the eqs.\mref{w6}:
\be\label{w8}
\mb[A,v_-]=\frac{i}{\kappa}(I-A)&&[A^*,v_-]=\frac{i}{\kappa}(A^*-{A^*}^2)\nn\\
\mb[A,v_+]=\frac{i}{\kappa}(A-A^2)&&[A^*,v_+]=\frac{i}{\kappa}(I-A^*)\nn\\
\mb[v_+,v_-]=\frac{i}{\kappa}(v_--v_+)&&[A,A^*]=0\nn\\
AA^*=A^*A=I\nn\\
\Delta A=A\otimes A&&\Delta A^*=A^*\otimes A^*\\
\Delta v_+=A\otimes v_++v_+\otimes I&&\Delta v_-=A^*\otimes v_-+v_-\otimes I\nn\\
S(A)=A^*&&S(A^*)=A\nn\\
S(v_+)=-A^*v_+&&S(v_-)=-Av_-\nn\\
\epsilon(A)=\epsilon(A^*)=1&&\epsilon(v_+)=\epsilon(v_-)=0\nn
\ee
\section{The left covariant 3D differential calculus\label{s3}}
In this section we describe 3D left covariant differential calculus on $E_\kappa(2)$.
This calculus results from the contraction of 3D left covariant differential
calculus on $SU_\mu(2)$, the latter being constructed by S.L.Woronowicz, \cite{b10}.

The right ideal $R$\ which defines the 3D calculus on $SU_\mu(2)$\ is generated by the following six elements, \cite{b10}:
$\rho^2;\;{\rho^*}^2;\;\hbox{$\rho^*\rho=\rho\rho^*$},\;(\sigma-I)(\rho-\rho^*);\;
(\sigma-I)(\rho+\rho^*);\;\sigma^*+\mu^*\sigma-(1+\mu^2)I$. Denoting by $\tilde{R_0}$\
the contraction of $R$\ we have:
\begin{moje}
\mb\\
The right ideal $\tilde{R_0}$\ is generating by the following elements:
\be\label{w9}
&&(a_0-I)(a_0^*-I);\;\;\;(a_0-I)w_0^*+\frac{1}{2\kappa}(a_0-a_0^*);\nn\\
&&(a_0-I)w_0-\frac{1}{2\kappa}(a_0-a_0^*);\;\;\;(a_0^*-I)w_0^*-\frac{1}{2\kappa}(a_0-a_0^*);\\
&&(a_0^*-I)w_o+\frac{1}{2\kappa}(a_0-a_0^*);\;\;\;w_0^*w_0+\frac{1}{2\kappa}(w_0-3w_0^*)-\frac{1}{4\kappa^2}(a_0-a_0^*);\nn\\
&&w_0^2+\frac{1}{2\kappa}(w_0^*-3w_0)+\frac{1}{4\kappa^2}(a_0-a_0^*);\;\;\;w_0^*-\frac{1}{2\kappa}(w_0-3w_0^*)+\frac{1}{4\kappa^2}(a_0-a_0^*)\nn
\ee
\end{moje}

The proof
is given in Appendix B.

It follows from the theorem  that the elements $a_0-a_0^*,\;\;w_0$\ and $w_0^*$\ from a~basis
in $\ker\epsilon/\tilde{R_0}$.
According to the Woronowicz theory, \cite{b15}, the space of left invariant 1--forms on
$\tilde{E}_\kappa(2)$\ is spanned by:
\be\label{w10}
\tilde{\varphi_0}&=&\pi r^{-1}(I\otimes(a_0-a_0^*))=a_0^*da_0-a_0da_0^*\nn\\
\tilde{\varphi_1}&=&\pi r^{-1}(I\otimes w_0)=-a_0^*w_0a_0da_0^*+a_0^*dw_0\\
\tilde{\varphi_2}&=&\pi r^{-1}(I\otimes w_0)=-a_0w_0^*a_0^*da_0+a_0dw_0^*\nn
\ee
Let us pass to $E_\kappa(2)$. The following corollary is a~simple consequence of
Theorem~1.
\newtheorem{moje2}{Corollary}
\begin{moje2}
\mb\\
The right ideal $R_0=E_\kappa(R)\cap \tilde{R_0}$\ is generated by the following
elements:
\be\label{w11}
&&(A-I)(A^*-I)\nn\\
&&(A-I)v_--\frac{i}{\kappa}(A^*-I);\;\;\;\;\;(A^*-I)v_-+\frac{i}{\kappa}(A^*-I);\nn\\
&&(A-I)v_++\frac{i}{\kappa}(A-I);\;\;\;\;\;(A^*-I)v_+-\frac{i}{\kappa}(A-I);\nn\\
&&v_-v_++\frac{i}{2\kappa}(v_++3v_-)+\frac{1}{4\kappa^2}(A^*-A)\\
&&v_+^2+\frac{i}{2\kappa}(3v_++v_-)+\frac{1}{4\kappa^2}(A^*-A)\nn\\
&&v_-^2+\frac{i}{2\kappa}(v_++3v_-)+\frac{1}{4\kappa^2}(A^*-A)\nn
\ee
\end{moje2}
Accordingly,  the elements: $A-A^*,\;\;v_+$\ and $v_-$\ form
a~basis in $\ker\epsilon/R_0$\ and the space of left invariant 1--forms is spanned
by the following forms:
\be\label{w12}
\varphi_0&=&\pi r^{-1}(I\otimes (A-A^*))=A^*dA-AdA^*\nn\\
\varphi_1&=&\pi r^{-1}(I\otimes v_+)=A^*dv_+\\
\varphi_2&=&\pi r^{-1}(I\otimes v_-)=Adv_-\nn
\ee
The following relations hold for the forms
defined by eqs.\mref{w10} and~\mref{w12}:
\be\label{w13}
\varphi_0&=&2\tilde{\varphi_0}\nn\\
\varphi_1&=&-\frac{i}{2\kappa}\tilde{\varphi_0}-i\tilde{\varphi_1}\\
\varphi_2&=&i\tilde{\varphi_2}\nn
\ee

Now one can easily find the following commutation rules between the invariant
forms and generators of $E_\kappa(2)$:
\be\label{w14}
\mb[A,\varphi_0]&=&0\;\;\;\;\;\;\;\;[A^*,\varphi_0]=0\nn\\
\mb[A,\varphi_1]&=&0\;\;\;\;\;\;\;\;[A^*,\varphi_1]=0\nn\\
\mb[A,\varphi_2]&=&0\;\;\;\;\;\;\;\;[A^*,\varphi_2]=0\nn\\
\mb[v_-,\varphi_0]&=&\frac{i}{\kappa}A^*\varphi_0\nn\\
\mb[v_-,\varphi_1]&=&\frac{i}{2\kappa}A^*(3\varphi_1+\varphi_2)-\frac{1}{4\kappa^2}A^*\varphi_0\nn\\
\mb[v_-,\varphi_2]&=&\frac{i}{2\kappa}A^*(\varphi_1+3\varphi_2)-\frac{1}{4\kappa^2}A^*\varphi_0\\
\mb[v_+,\varphi_0]&=&\frac{i}{\kappa}A\varphi_0\nn\\
\mb[v_+,\varphi_1]&=&\frac{i}{2\kappa}A(3\varphi_1+\varphi_2)-\frac{1}{4\kappa^2}A\varphi_0\nn\\
\mb[v_+,\varphi_2]&=&\frac{i}{2\kappa}A(\varphi_1+3\varphi_2)-\frac{1}{4\kappa^2}A\varphi_0\nn
\ee

To complete the description of the first order differential calculus we must introduce
the $*$--operator. It is easy to see that the involution acts as follow:
\be\label{w15}
\varphi_0^*=-\varphi_0,&\varphi_1^*=\varphi_2,&\varphi_2^*=\varphi_1
\ee

The next step is the construction of the higher order differential calculus.
Because our first order calculus is not bicovariant the general Woronowicz theory
does not work in this case. However we can apply our contraction procedure
to obtain the higher differential forms from the ones on $SU_\mu(2)$\ constructed
by Woronowicz, \cite{b10}. Let us recall that the left
invariant basic forms on $SU_\mu(2)$\ introduced in Ref.\cite{b10} read:
\be\label{w16}
\Omega_0&=&\rho^*d\sigma^*-\sigma d\rho^*\nn\\
\Omega_1&=&\sigma^*d\sigma-\rho^*d\rho\\
\Omega_2&=&\rho d\sigma-\mu^{-1}\sigma d\rho\nn
\ee
One can easily find the following expansions:
\be\label{w17}
\Omega_0&=&\frac{1}{4}\varphi_0+\frac{1}{2R}\ll[i(\varphi_1-\varphi_2)-\frac{1}{\kappa}\varphi_0\r]+O\ll(\frac{1}{R^2}\r)\nn\\
\Omega_1&=&\frac{1}{2R}\ll[-i(\varphi_1+\varphi_2)+\frac{1}{2\kappa}\varphi_0\r]+O\ll(\frac{1}{R^2}\r)\\
\Omega_2&=&-\frac{1}{4}\varphi_0+\frac{1}{2R}\ll[i(\varphi_1-\varphi_2)+\frac{5}{4\kappa}\varphi_0\r]+O\ll(\frac{1}{R^2}\r)\nn
\ee
Inserting these expansions into the following external product identies written
out in~\cite{b10}:
\be\label{w18}
\Omega_0\wedge\Omega_0&=&0\;\;\;\;\;\;\;\;\;\Omega_1\wedge\Omega_1=0\nn\\
\Omega_1\wedge\Omega_0&=&-\mu^4\Omega_0\wedge\Omega_1\\
(\Omega_2+\Omega_0)\wedge\Omega_0&=&-\mu^2\Omega_0\wedge(\Omega_0+\Omega_2)\nn\\
(\Omega_0+\Omega_2)\wedge\Omega_1&=&-\mu^{-4}\Omega_1\wedge\Omega_0-\mu^4\Omega_1\wedge\Omega_2\nn
\ee
we obtain the following external product rules:
\be\label{w19}
&&\varphi_0\wedge\varphi_0=0\nn\\
&&\varphi_0\wedge\varphi_1+\varphi_1\wedge\varphi_0=0\nn\\
&&\varphi_0\wedge\varphi_0+\varphi_2\wedge\varphi_0=0\\
&&\varphi_2\wedge\varphi_1+\varphi_1\wedge\varphi_2+\frac{i}{4\kappa}\varphi_2\wedge\varphi_0-\frac{i}{4\kappa}\varphi_1\wedge\varphi_0=0\nn\\
&&\varphi_2\wedge\varphi_2-\frac{3i}{8\kappa}\varphi_1\wedge\varphi_0-\frac{5i}{8\kappa}\varphi_2\wedge\varphi_0=0\nn\\
&&\varphi_1\wedge\varphi_1+\frac{5i}{8\kappa}\varphi_1\wedge\varphi_0+\frac{3i}{8\kappa}\varphi_2\wedge\varphi_0=0\nn
\ee

The following Cartan--Maurer equations complete the description of the differential calculus:
\be\label{w20}
d\varphi_0&=&0\nn\\
d\varphi_1&=&-\fr\varphi_0\wedge\varphi_1\\
d\varphi_2&=&\fr\varphi_0\wedge\varphi_2\nn
\ee

Woronowicz theory, \cite{b15}, provides us with the general construction of Lie
algebra once the left--covariant calculus is given (see also \cite{b10}). Following
general framework we introduce the counterparts of left invariant vector fields
by the formula:
\be\label{w21}
dx&=&(\chi_0*x)\varphi_0+(\chi_1*x)\varphi_1+(\chi_2*x)\varphi_2
\ee
where $x\in E_\kappa(2)$\ and $\chi*x=(id\otimes\chi)\Delta x$.

Applying  the exterior derivative to both side of eq.\mref{w21} and taking into account
Cartan--Maurer equations as well as and the exterior algebra relations~\mref{w19} we
arrive at the following commutation rules:
\be\label{w22}
\mb[\chi_1,\chi_0]&=&\frac{5i}{8\kappa}\chi_1^2-\fr\chi_1-\frac{i}{2\kappa}\chi_1\chi_2-\frac{3i}{8\kappa}\chi_2^2\nn\\
\mb[\chi_2,\chi_0]&=&\frac{3i}{8\kappa}\chi_1^2+\frac{i}{4\kappa}\chi_1\chi_2+\fr\chi_2-\frac{5i}{8\kappa}\chi_2^2\\
\mb[\chi_1,\chi_2]&=&0\nn
\ee

It is easy to check that the involution  acts as follow:
\be\label{w23}
\chi_0^*=\chi_0,&\chi_1^*=-\chi_2,&\chi_2^*=-\chi_1
\ee

The coproduct for the functionals $\chi_i,\; i=0,1,2$\ is defined by:
\be\label{w24}
\Delta\chi_i&=&\sum\limits^2_{j=0}\chi_j\otimes f_{ij}+I\otimes\chi_i
\ee
where the functionals $f_{ij}$\ are given by the relations:
\be\label{w25}
\varphi_i x&=&\sum\limits^2_{j=0}(f_{ij}*x)\varphi_j
\ee

On the other hand it has been shown (\cite{b11},\cite{b12}) that the deformed
Lie algebra $e_\kappa(2)$\ dual to $E_\kappa(2)$\ is a Hopf algebra generated
by three selfadjoint elements $P_1,P_2$\ and $J$\ subject to the following
relations:
\be\label{w26}
\mb[P_1,P_2]&=&0\nn\\
\mb[J,P_1]&=&=iP_2\nn\\
\mb[J,P_2]&=&-i\kappa\sinh\ll(\frac{P_1}{\kappa}\r)\\
\Delta P_1&=&I\otimes P_1+P_1\otimes I\nn\\
\Delta P_2&=&e^{-P_1/2\kappa}\otimes P_2+P_2\otimes e^{P_1/2\kappa}\nn\\
\Delta J&=&e^{-P_1/2\kappa}\otimes J+J\otimes e^{P_1/2\kappa}\nn
\ee

Now, one can pose the question what is the relation between the functionals $\chi_i,\;i=0,1,2$\
and the generators $P_1,P_2$\ and $J$\ of the $e_\kappa(2)$. It is not difficult to check that the answer to this
question is given by the following relations:
\be\label{w27}
\chi_0&=&\fr e^{P_1/2\kappa}(J+\frac{i}{2\kappa}P_2)-\frac{1}{8}(e^{2P_1/\kappa}-I)\nn\\
\chi_1&=&\frac{i\kappa}{4}(e^{2P_1/\kappa}-I)-\fr P_2e^{P_1/2\kappa}\nn\\
\chi_2&=&\frac{i\kappa}{4}(e^{2P_1/\kappa}-I)+\fr P_2e^{P_1/2\kappa}\nn\\
f_{01}&=&f_{02}=0\;\;\;\;\;\;\;\;\;\;f_{00}=e^{P_1/\kappa}\\
f_{20}&=&f_{10}=\frac{i}{4\kappa}(e^{2P_1/\kappa}-e^{-P_1/\kappa})\nn\\
f_{22}&=&f_{11}=\fr (e^{2P_1/\kappa}+e^{P_1/\kappa})\nn\\
f_{21}&=&f_{12}=\fr (e^{2P_1/\kappa}-e^{P_1/\kappa})\nn
\ee

\section{The bicovariant 4D$_+$\ differential calculus on $\tilde{E_\kappa}(2)$\label{s4}}
Let us recall, \cite{b13}, that the right ideal $R_+$\ (resp.$R_-$) which
defines the 4D$_+$\ (resp.4D$_-$) bicovariant differential calculus on
$SU_\mu(2)$\ is generated by the following elements:
\be\label{w28}
&&\rho^2;\;\;\;\;\;\rho(\sigma-\sigma^*);\;\;\;\;\;\sigma^2+\mu^2{\sigma^*}^2+(1+\mu^2)^2\rho\rho^*-(1+\mu^2)I;\nn\\
&&{\rho^*}^2;\;\;\;\;\;\rho^*(\sigma-\sigma^*),\;a\rho;\;\;\;\;\;a\rho^*;\;\;\;\;\;a(\sigma-\sigma^*);\\
&&a(\mu^2\sigma+\sigma^*-(1+\mu^2)I)\nn
\ee
where $a=\mu^2\sigma+\sigma^*-[(1+\mu^4)/\mu]I$\
(resp.$a=\mu^2\sigma+\sigma^*+[(1+\mu^4)/\mu]I$).

Let us denote by $\tilde{R}_{0+}$\ the contraction of the ideal
$R_+$. then we have:
\begin{moje}
\mb\\
The right ideal $\tilde{R_{0+}}$\ is generated by the following
elements:
\be
a_0+a_0^*-2;&(a_0-I)w_0;&(a_0^*-I)w_0\nn\\
(a_0-I)w_0^*;&(a_0^*-I)w_0^*;&w_0^2+\frac{1}{\kappa}w_0;\nn\\
w_0^*-\frac{1}{\kappa}w_0^*\nn
\ee
\end{moje}

For the proof
see Appendix C.

The  ideal $\tilde{R}_{0+}$\ is ad--invariant then so we
can apply the general Woronowicz theory, \cite{b15}.
It
follows from the theorem~3 that the space of left invariant
1--forms is spanned by the following 1--forms:
\be\label{w29}
\psi_1&=&\pi r^{-1}(I\otimes(a_0-a_0^*))=a_0^*da_0-a_0da_0^*=2a_0^*da_0=-2a_0da_0^*\nn\\
\psi_2&=&\pi r^{-1}(I\otimes w_0)=\fr a_0^*w_0\psi_1+a_0^*dw_0\nn\\
\psi_3&=&\pi r^{-1}(I\otimes w_0^*)=-\fr a_0w_0^*\psi_1+a_0dw_0^*\\
\psi_4&=&\pi r^{-1}(I\otimes w_0w_0^*)=-w_0a_0^*\psi_3-\frac{1}{2\kappa}(1-a_0^2)(\psi_2-\frac{1}{2\kappa}\psi_1)+d(w_0w_0^*)\nn
\ee

One can easily find the following commutation relations:
\be\label{w30}
\mb[a_0,\psi_1]&=&0\;\;\;\;\;\;\;\;\;[w_0,\psi_1]=0\nn\\
\mb[a_0^*,\psi_1]&=&0\;\;\;\;\;\;\;\;\;[w_0^*,\psi_1]=0\nn\\
\mb[a_0,\psi_2]&=&\frac{1}{2\kappa}a_0\psi_1\nn\\
\mb[a_0^*,\psi_2]&=&-\frac{1}{2\kappa}a_0^*\psi_1\nn\\
\mb[w_0,\psi_2]&=&-\frac{1}{2\kappa}w_0\psi_1+\frac{1}{\kappa}a_0\psi_2\nn\\
\mb[w_0^*,\psi_2]&=&\frac{1}{2\kappa}w_0^*\psi_1-a_0^*\psi_4\nn\\
\mb[a_0,\psi_3]&=&-\frac{1}{2\kappa}a_0\psi_1\\
\mb[a_0^*,\psi_3]&=&\frac{1}{2\kappa}a_0^*\psi_1\nn\\
\mb[w_0,\psi_3]&=&\frac{1}{2\kappa}w_0\psi_1-a_0\psi_4-\frac{1}{\kappa}a_0(\psi_1+\psi_3)\nn\\
\mb[w_0^*,\psi_3]&=&-\frac{1}{2\kappa}w_0^*\psi_1-\frac{1}{\kappa}a_0^*\psi_3\nn\\
\mb[a_0,\psi_4]&=&\frac{1}{2\kappa^2}a_0\psi_1\nn\\
\mb[a_0^*,\psi_4]&=&-\frac{1}{2\kappa^2}a_0^*\psi_1\nn\\
\mb[w_0,\psi_4]&=&-\frac{1}{2\kappa^2}w_0\psi_1+\frac{1}{\kappa^2}a_0\psi_2\nn\\
\mb[w_0^*,\psi_4]&=&\frac{1}{2\kappa^2}w_0^*\psi_1+\frac{1}{\kappa}a_0^*\psi_4+\frac{2}{\kappa^2}a_0^*\psi_3\nn
\ee

The external product identies read:
\be\label{w31}
&&\psi_1\wedge\psi_1=0\nn\\
&&\psi_1\wedge\psi_2+\psi_2\wedge\psi_1=0\nn\\
&&\psi_1\wedge\psi_3+\psi_3\wedge\psi_1=0\nn\\
&&\psi_1\wedge\psi_4+\psi_4\wedge\psi_1=0\nn\\
&&\psi_2\wedge\psi_2-\frac{1}{\kappa}\psi_1\wedge\psi_2=0\\
&&\psi_2\wedge\psi_3+\psi_3\wedge\psi_2+\frac{1}{\kappa}\psi_1\wedge(\psi_3+\psi_2)=0\nn\\
&&\psi_2\wedge\psi_4+\psi_4\wedge\psi_2=0\nn\\
&&\psi_3\wedge\psi_3-\frac{1}{\kappa}\psi_1\wedge\psi_3=0\nn\\
&&\psi_3\wedge\psi_4+\psi_4\wedge\psi_3+\frac{1}{\kappa^2}\psi_1\wedge(\psi_3-\psi_2)=0\nn\\
&&\psi_4\wedge\psi_4+\frac{1}{\kappa^3}\psi_1\wedge\psi_2=0\nn
\ee

while Cartan--Maurer equations are given by:
\be\label{w32}
d\psi_1&=&0\nn\\
d\psi_2&=&-\psi_1\wedge\psi_2\\
d\psi_3&=&\psi_1\wedge\psi_3\nn\\
d\psi_4&=&\frac{1}{\kappa}\psi_1\wedge\psi_2\nn
\ee

The quantum Lie algebra reads(
$dx=(\xi_1*x)\psi_1+(\xi_2*x)\psi_2+(\xi_3*x)\psi_3+(\xi_4*x)\psi_4\nn$):

\be\label{w33}
\mb[\xi_1,\xi_2]&=&-\frac{1}{\kappa}\xi_2^2+\frac{1}{\kappa}\xi_3\xi_2+\xi_2-\frac{1}{\kappa^2}\xi_4\xi_3+\frac{1}{\kappa^3}\xi_4^2-\frac{1}{\kappa}\xi_4\nn\\
\mb[\xi_1,\xi_3]&=&\frac{1}{\kappa}\xi_3\xi_2-\frac{1}{\kappa}\xi_3^2+\frac{1}{\kappa^2}\xi_4\xi_3-\xi_3\nn\\
\mb[\xi_1,\xi_4]&=&0\\
\mb[\xi_2,\xi_3]&=&0\nn\\
\mb[\xi_2,\xi_4]&=&0\nn\\
\mb[\xi_3,\xi_4]&=&0\nn
\ee

\section{4D$_-$\ differential calculus on $\tilde{E}_\kappa(2)$\label{s5}}

Denoting by $\tilde{R}_{0-}$\ the contraction of the ideal $R_-$\ (see eqs.\mref{w28}) we have:
\begin{moje}
\mb\\
The right ideal $\tilde{R}_{0-}$\ is generated by the following elements:
\be
a_0^2+a_0-a_0^*-I;&&{a_0^*}^2+a_0^*-a_0-I;\nn\\
(a_0+I)w_0;&&(a_0^*+I)w_0;\nn\\
(a_0+I)w_0^*;&&(a_0^*+I)w_0^*;\nn\\
w_0^2-\frac{1}{\kappa}w_0;&&{w_0^*}^2+\frac{1}{\kappa}w_0^*\nn\\
&&w_0w_0^*-\frac{1}{\kappa}w_0-\frac{1}{\kappa^2}(a_0^*-I)\nn
\ee
\end{moje}
Proof:\hfill\mb\\
Because the proof of this theorem is very similar to the proof of theorem~3 and will be omitted.

The ideal $\tilde{R}_{0-}$\ is ad--invariant and we can follow the Woronowicz theory. The space of
left invariant 1--forms is spanned by the following four 1--forms:
\be\label{w34}
\Phi_1&=&\pi r^{-1}(I\otimes(a_0-I))=a_0^*da_0\nn\\
\Phi_2&=&\pi r^{-1}(I\otimes(a_0^*-I))=a_0da_0^*\\
\Phi_3&=&\pi r^{-1}(I\otimes w_0)=a_0^*dw_0-a_0^*w_0a_0da_0^*\nn\\
\Phi_4&=&\pi r^{-1}(I\otimes w_0^*)=a_0dw_0^*-a_0w_0^*a_0^*da_0\nn
\ee

The following commutation relations hold between above forms and the generators of
$\tilde{E}_\kappa(2)$:
\be\label{w35}
\Phi_1a_0&=&a_0(\Phi_2-2\Phi_1)\nn\\
\Phi_1a_0^*&=&-a_0^*\Phi_2\nn\\
\Phi_2a_0&=&-a_0\Phi_1\nn\\
\Phi_2a_0^*&=&a_0^*(\Phi_1-2\Phi_2)\nn\\
\Phi_3a_0&=&-a_0\Phi_3+\frac{1}{2\kappa}a_0(\Phi_1-\Phi_2)\nn\\
\Phi_3a_0^*&=&-a_0^*\Phi_3+\frac{1}{2\kappa}a_0^*(\Phi_2-\Phi_1)\nn\\
\Phi_4a_0^*&=&-a_0\Phi_4+\frac{1}{2\kappa}a_0(\Phi_2-\Phi_1)\nn\\
\Phi_4a_0^*&=&-a_0^*\Phi_4+\frac{1}{2\kappa}a_0^*(\Phi_1-\Phi_2)\\
\Phi_1w_0&=&-w_0\Phi_2-2a_0\Phi_3\nn\\
\Phi_1w_0^*&=&w_0^*\Phi_2-2w_0^*\Phi_1-2a_0^*\Phi_4\nn\\
\Phi_2w_0&=&w_0\Phi_1-2w_0\Phi_2-2a_0\Phi_3\nn\\
\Phi_2w_0^*&=&-w_0^*\Phi_1-2a_0^*\Phi_4\nn\\
\Phi_3w_0&=&-w_0\Phi_3-\frac{1}{2\kappa}w_0(\Phi_1-\Phi_2)+\frac{1}{\kappa}a_0\Phi_3\nn\\
\Phi_3w_0^*&=&-w_0^*\Phi_3+\frac{1}{2\kappa}w_0^*(\Phi_1-\Phi_2)+\frac{1}{\kappa}a_0^*\Phi_3+\frac{1}{\kappa^2}a_0^*\Phi_2\nn\\
\Phi_4w_0&=&-w_0\Phi_4+\frac{1}{2\kappa}w_0(\Phi_1-\Phi_2)-\frac{1}{\kappa}a_0\Phi_4+\frac{1}{\kappa^2}a_0\Phi_2\nn\\
\Phi_4w_0^*&=&-w_0^*\Phi_4+\frac{1}{2\kappa}w_0^*(\Phi_2-\Phi_1)-\frac{1}{\kappa}a_0^*\Phi_4\nn
\ee

The external product identies read:
\be\label{w36}
&&\Phi_1\wedge\Phi_1=0\nn\\
&&\Phi_2\wedge\Phi_2=0\nn\\
&&\Phi_1\wedge\Phi_2+\Phi_2\wedge\Phi_1=0\nn\\
&&\Phi_3\wedge\Phi_1+3\Phi_1\wedge\Phi_3-2\Phi_2\wedge\Phi_3=0\nn\\
&&\Phi_3\wedge\Phi_2-\Phi_2\wedge\Phi_3+2\Phi_1\wedge\Phi_3\\
&&\Phi_3\wedge\Phi_3-\frac{1}{\kappa}\Phi_1\wedge\Phi_3+\frac{1}{\kappa}\Phi_2\wedge\Phi+3=0\nn\\
&&\Phi_4\wedge\Phi_1-\Phi_1\wedge\Phi_4+2\Phi_2\wedge\Phi_4=0\nn\\
&&\Phi_4\wedge\Phi2+3\Phi_2\wedge\Phi_4-2\Phi_1\wedge\Phi_4=0\nn\\
&&\Phi_4\wedge\Phi_3+\Phi_3\wedge\Phi_4+\frac{1}{\kappa}(\Phi_1-\Phi_2)\wedge\Phi_3+\frac{1}{\kappa}(\Phi_1-\Phi_2)\wedge\Phi_4=0\nn\\
&&\Phi_4\wedge\Phi_4-\frac{1}{\kappa}(\Phi_1-\Phi_2)\wedge\Phi_4=0\nn
\ee

The following Cartan--Maurer formulas complete the description of the second order
calculus:
\be\label{w37}
d\Phi_1&=&0\nn\\
d\Phi_2&=&0\\
d\Phi_3&=&(\Phi_1-\Phi_2)\wedge\Phi_3\nn\\
d\Phi_4&=&(\Phi_2-\Phi_1)\wedge\Phi_4\nn
\ee

The Lie algebra relations
($dx=(\eta_1*x)\Phi_1+(\eta_2*x)\Phi_2+(\eta_3*x)\Phi_3+(\eta_4*x)\Phi_4$) read:
\be\label{w38}
\mb[\eta_1,\eta_2]&=&0\nn\\
\mb[\eta_1\eta_3]&=&2\eta_3\eta_1+2\eta_3\eta_2-\frac{1}{\kappa}\eta_3^2+\frac{1}{\kappa}\eta_4\eta_3-\eta_3\nn\\
\mb[\eta_1,\eta_4]&=&-2\eta_4\eta_1-2\eta_4\eta_2+\frac{1}{\kappa}\eta_4\eta_3-\frac{1}{\kappa}\eta_4^2+\eta_4\\
\mb[\eta_2,\eta_3]&=&-2\eta_3\eta_1-2\eta_3\eta_2+\frac{1}{\kappa}\eta_3^2-\frac{1}{\kappa}\eta_4\eta_3+\eta_3\nn\\
\mb[\eta_2,\eta_4]&=&2\eta_4\eta_1+2\eta_4\eta_2-\frac{1}{\kappa}\eta_4\eta_3+\frac{1}{\kappa}\eta_4^2-\eta_4\nn\\
\mb[\eta_3,\eta_4]&=&0\nn
\ee

\section{The bicovariant 4D differential calculus on $E_\kappa(2)$\label{s6}}

It is known, \cite{b14}, that there exists only one four dimensional bicovariant
calculus on $E_\kappa(2)$. On the other hand, as it was shown in sec.\ref{s4} and~\ref{s5} 
there exist two fourdimensional calculi on $\tilde{E}_\kappa(2)$, so they
have to correspond to the same calculus on $E_\kappa(2)$. Indeed, we have:
\begin{moje}
\mb\\
The following equalities hold:
\be
\tilde{R}_{0+}\cap E_\kappa(2)=&\tilde{R}_{0-}\cap E_\kappa(2)&=R\nn
\ee
where the right ideal $R$\ defines the 4D bicovariant differential calculus
described in~\cite{b14}.
\end{moje}
Proof.\\
The proof of this theorem is straightforward.
\vskip1cm

Let us recall the basic left invariant 1--forms, introduced in Ref.\cite{b14}.
\be\label{w39}
w_0&=&\fr(A^*dA-AdA^*)=A^*dA=-AdA^*\nn\\
w_+&=&A^*dv_+\\
w_-&=&A^*dv_-\nn\\
\tilde{w}_0&=&d(v_+v_-+\frac{i}{\kappa}v_+)-v_+dv_--v_-dv_+\nn
\ee

Like in undeformed case, we can express the left invariant 1--forms and the left
invariant fields defined on the $E_\kappa(2)$\ by the ones defined on the 
$\tilde{E}_\kappa(2)$. For D$_+$\ we get ($dx=(\chi_0*x)\omega_0+(\chi_+*x)\omega_++(\chi_-*x)\omega_-)$.

\be\label{w40}
w_0&=&\psi_1\nn\\
w_+&=&-i\psi_2\nn\\
w_-&=&i\psi_3-\frac{i}{2\kappa}\psi_1\nn\\
\tilde{w}_0&=&\psi_4+\frac{1}{\kappa}\psi_2+\frac{1}{2\kappa^2}\psi_1\\
\chi_0&=&\xi_1\nn\\
\chi_+&=&i(\xi_2-\frac{1}{\kappa}\xi_4)\nn\\
\chi_-&=&-i\xi_3\nn\\
\tilde{\chi}_0&=&\xi_4\nn
\ee
while for D$_-$\ we obtain:
\be\label{w41}
w_0&=&\Phi_2-\Phi_1\nn\\
w_+&=&i\Phi_3\nn\\
w_-&=&-i\Phi_4+\frac{i}{\kappa}(\Phi_1-\Phi_2)\nn\\
\tilde{w}_0&=&\frac{3}{2\kappa^2}\Phi_2-\frac{1}{2\kappa^2}\Phi_1\\
\chi_0&=&\fr (\eta_2-\eta_1)\nn\\
\chi_+&=&-i\eta_3\nn\\
\chi_-&=&i\eta_4\nn\\
\tilde{\chi}_0&=&\kappa^2(\eta_2+\eta_1)\nn
\ee

Let us conclude with the following remark. In the classical case the differential
calculus is obtained with the choice $R=(\ker\epsilon)^2$, i.e. it is determined
by the ideal consisting of functions that vanish, up to second order, at the group identity.
Therefore, local diffeomorphism gives unique relation between differential calculi.
However, in the quantum case situation looks differently. Two different calculi
on $\tilde{E}_\kappa(2)$\ reduce to the single one on $E_\kappa(2)$. In order
to get some insight let us consider Hopf subalgebra of $\tilde{E}_\kappa(2)$\ 
generated by $a_0,a_0^*$. It is commutative Hopf algebra so we can speak
in terms of algebra of functions on $U(1)$. In the D$_+$\ case we obtain
the standard calculus on $U(1)$. Indeed, denoting $a_0=a^{i\Theta}$\ we see
that the corresponding ideal is generated by $\cos\Theta-1\approx\Theta^2$.
On the other hand in the D$_-$\ case the ideal is generated by $\cos 2\Theta-1$\ 
and $\sin\Theta(\cos\Theta+1)$. Therefore, it consists of functions vanishing not 
only at $\Theta=0$\ but also at $\Theta=\Pi$. However, under the mapping $a_0\to A=a_0^2$,
which is double covering, it procedures the ideal of functions vanishing at the group
identity (in a special way). We see that in the quantum case, generically, the relation
between calculi depends on global properties of the mapping.
\section{Appendix}

{\bf A)} Proof of the theorem 1.\\
We begin from the analysis of the relations
given in eq.\mref{w1} (up to terms of order $\frac{1}{R}$).\\

$\rho^*\rho=\rho^*\rho$\ implies:
\be\label{w42}
\lim\limits_{R\to\infty}R^n[w+w^*,a-a^*]=0&&n=0,1,2,\ldots
\ee

$\mu(\rho-\rho^*)\sigma=\sigma(\rho-\rho^*)$\ implies:
\be
&&\mb[a_0,a_0^*]=0\label{w43}\\
&&\lim\limits_{R\to\infty}(2R[a,a^*]+[a-a^*,w-w^*]+\frac{1}{\kappa}(a-a^*)(a+a^*))=0\label{w44}
\ee

$\mu(\rho+\rho^*)\sigma=\sigma(\rho+\rho^*)$\ gives:
\be
&&\mb[w_0+w_0^*,a_0+a_0^*]=0\label{w45}\\
&&\lim\limits_{R\to\infty}(R[w+w^*,a+a^*]+\frac{1}{\kappa}[w+w^*,a]+\nn\\
&&\mb+\frac{1}{\kappa}(w+w^*)(a+a^*)+2[w,w^*])=0\label{w46}
\ee

$\sigma^*\sigma+\rho\rho^*=I$\ gives:
\be
&&a_0^*a_0+a_0a_0^*-2=0\label{w47}\\
&&\lim\limits_{R\to\infty}(2R(a^*a+a^*a-2)+\frac{1}{\kappa}(a+a^*)^2+[a+a^*,w^*-w]+\nn\\
&&\mb+[a-a^*,w+w^*])=\frac{4}{\kappa}\label{w48}
\ee

Finally, from the relation $\sigma\sigma^*+\mu^2\rho^*\rho=I$\ we obtain:
\be\label{w49}
&&\lim\limits_{R\to\infty}(2R(aa^*+a^*a-2)+\frac{1}{\kappa}(a+a^*)^2+[a+a^*,w-w^*]-\nn\\
&&\mb-\frac{2}{\kappa}(a-a^*)^2+[a-a^*,w+w^*])=\frac{4}{\kappa}
\ee

From eqs.\mref{w43},\mref{w47},\mref{w42} and~\mref{w45} one concludes:
\be\label{w50}
&a_0a+0^*=a_0^*a_0=I&
\ee
and
\be\label{w51}
\mb[w_0,a_0]&=&-[w_0^*,a_0]\\
\mb[w_0,a_0^*]&=&[w_0^*,a_0^*]\nn
\ee

If we subtract eq.\mref{w49} from eq.\mref{w48} and take into account eq.\mref{w51}
we arrive at:
\be\label{w52}
\mb[a_0,w_0]+[a_0^*,w_0]&=&\frac{1}{2\kappa}(a_0^@-{a_0^*}^2-a_0^2)
\ee

Subtracting from eq.\mref{w44} its conjugation and using eq.\mref{w51} we get:
\be\label{w53}
-[a_0,w_0]+[a_0^*,w_0]&=&\frac{1}{2\kappa}({a_0^*}^2-a_0^2)
\ee

Now adding to eq.\mref{w46} its conjugation and again using eq.\mref{w51} we arrive
at the relation:
\be\label{w54}
\mb[w_0,w_0^*]&=&-\frac{1}{2\kappa}(a_0+a_0^*)(w_0+w_0^*)
\ee

Finally from eqs.\mref{w50},\mref{w52},\mref{w53} and~\mref{w54} we get the
commutation rules~\mref{w6}.

{\bf B)} The proof of the theorem 2.\\
From the inclusions $\rho^2\in R,\;{\rho^*}^2\in R$\ and $\rho\rho^*\in R$\ 
we immediately obtain that
\be\label{w55}
(a_0-a_0^*)^2\in\tilde{R}_0;\;(w_0+w_0^*)^2\in\tilde{R}_0\;\mbox{and}\;(a_0-a_0^*)(w_0+w_0^*)\in\tilde{R}_0
\ee

Due to $(\sigma-I)(\rho\pm\rho^*)\in R$, 
$\sigma^*+\mu^2\sigma-(1+\mu^2)I\in R$\ and  the relations
\mref{w55} one gets:
\be
&&a_0+a_0^*-2\in \tilde{R}_0\label{w56}\\
&&\lim\limits_{R\to\infty}R(a+a^*-2)\in\tilde{R}_0\label{w57}\\
&&\lim\limits_{R\to\infty}(R^2(a+a^*-2)-\frac{1}{\kappa^2}(a-I)+\frac{1}{\kappa}(w^*-3w))\in\tilde{R}_0\label{w58}\\
&&w_0^*(w_0+w_0^*)\in\tilde{R}_0\label{w59}\\
&&w_0(w_0+w_0^*)\in\tilde{R}_0\label{w60}\\
&&(a_0-I)w_0^*+\frac{1}{2\kappa}(a_0-a_0^*)\in\tilde{R}_0\label{w61}\\
&&(a_0-I)w_0-\frac{1}{2\kappa}(a_0-a_0^*)\in\tilde{R}_0\label{w62}\\
&&(a_0^*-I)w_0^*-\frac{1}{2\kappa}(a_0-a_0^*)\in\tilde{R}_0\label{w63}\\
&&(a_0^*-I)w_0+\frac{1}{2\kappa}(a_0-a_0^*)\in\tilde{R}_0\label{w64}
\ee

From the relations $\sigma^*\sigma+\rho^*\rho=I,\;\sigma\sigma^*+\mu^2\rho^*\rho=I$\ and
the inclusion $\rho\rho^*\in R$\ it follows that $[\sigma,\sigma^*]\in R$\ and
$\sigma\sigma^*-I\in R$.

The inclusion $[\sigma,\sigma^*]\in R$\ gives 
\be\label{w65}
\lim\limits_{R\to\infty}(\frac{1}{\kappa}R[a,a^*]+R[w^*-w,a+a^*]+\frac{1}{\kappa}[w-w^*,a+a^*])\in\tilde{R}_0
\ee

Taking into account eq.\mref{w57} we obtain from  $\sigma\sigma^*-I\in R$\ the relations:
\be\label{w66}
&&\lim\limits_{R\to\infty}(R^2(a+a^*-2)(a+a^*+2)+\frac{R}{\kappa}[a,a^*]+\nn\\
&&+R[w^*-w,a+a^*]+\frac{1}{\kappa}[a,w]+\frac{1}{\kappa}[w^*,a^*]-\\
&&-\frac{2}{\kappa}aw^*-\frac{2}{\kappa}wa^*-(w^*-w)^2)\in\tilde{R}_0\nn
\ee

Eqs.\mref{w55},\mref{w56},\mref{w58}--\mref{w61} and~\mref{w64}-\mref{w66} imply:
\be\label{w67}
w_0^*w_0+\frac{1}{2\kappa}(w_0-3w_0^*)-\frac{1}{4\kappa^2}(a_0-a_0^*)\in\tilde{R}_0
\ee

Finally, using again eqs.\mref{w59} and~\mref{w60} we get
\be
&&w_0^2+\frac{1}{2\kappa}(w_0^*-3w_0)+\frac{1}{4\kappa^2}(a_0-a_0^*)\in\tilde{R}_0\label{w68}\\
&&{w_0^*}^2-\frac{1}{2\kappa}(w_0-3w_0^*)+\frac{1}{4\kappa^2}(a_0-a_0^*)\in\tilde{R}_0\label{w69}
\ee

The relations \mref{w56},\mref{w61}-\mref{w64} and~\mref{w67}--\mref{w69} completely describe th ideal
$\tilde{R}_0$.

{\bf C)} Proof of the theorem 3.\\
The proof of this theorem is similar to the proof of theorem~2. After long and tedious
analysis of the generators of the ideal $R_+$\ we conclude  that:
\be
&&a_0+a_0^*-2\in\tilde{R}_{0+}\nn\\
&&(a_0-I)w_0\in\tilde{R}_{0+}\;\;\;(a_0-I)w_0^*\in\tilde{R}_{0+}\label{w70}\\
&&(a_0^*-I)w_0\in \tilde{R}_{0+}\;\;\;(a_0^*-I)w_0^*\in\tilde{R}_{0+}\nn\\
&&\mb\nn\\
&&\lim\limits_{R\to\infty}(2R^2(aa^*+a^*a-2)+R(\frac{1}{\kappa}(a^2+{a^*}^2+2a^*a-4)+\nn\\
&&+[a+a^*,w^*-w])+\frac{1}{2\kappa^2}(a^2+{a^*}^2+4a^*a-6)+\frac{1}{\kappa}([w,a]+\label{w71}\\
&&+[a^*,w^*]-2w^*a-2a^*w)-(w^*-w)^2+(w+w^*)^2)\in\tilde{R}_{0+}\nn\\
&&\mb\nn\\
&&\lim\limits_{R\to\infty}(-2R^2(a-a^*)^2-4R^2(2-aa^*-a^*a)-\frac{4}{\kappa^2}(a-a^*)^2+\nn\\
&&+\frac{1}{\kappa^2}(-4a^2+13aa^*+13a^*a-22)+2(w-w^*)^2+4(w+w^*)+\label{w72}\\
&&+\frac{1}{\kappa}(2\{a^*,w\}-4\{a^*,w^*\}-2\{a,w^*\}))\in\tilde{R}_{0+}\nn\\
&&\mb\nn\\
&&(w_0+w_0^*)(w_0^*-w_0)\in\tilde{R}_{0+}\label{w73}\\
&&\mbox{and}\nn\\
&&\lim\limits_{R\to\infty}R(a-a^*)^2\in\tilde{R}_{0+}\label{w74}
\ee

The relation $(\mu^2-1)(\rho+\rho^*)^2-2[\sigma,\sigma^*]\in R_+$\ gives:
\be\label{w75}
\lim\limits_{R\to\infty}R(\frac{1}{\kappa}[a,a^*]+[a+a^*,w-w^*])+\frac{1}{\kappa}[a_0+a_0^*,w_0-w_0^*]\in\tilde{R}_{0+}
\ee

while from eqs.\mref{w71},\mref{w72},\mref{w74} and~\mref{w75} it follows that:
\be\label{w76}
w_0^2+\frac{1}{\kappa}w_0+{w_0^*}^2-\frac{1}{\kappa}w_0\in\tilde{R}_{0+}
\ee

while from eqs.\mref{w76} and~\mref{w73} we obtain
\be
w_0^2+\frac{1}{\kappa}w_0\in\tilde{R}_{0+}\label{w77}\\
{w_0^*}^2-\frac{1}{\kappa}w_0^*\in\tilde{R}_{0+}\label{w78}
\ee

Eqs.\mref{w70},\mref{w77} and~\mref{w78} completely describe the ideal $\tilde{R}_{0+}$.


\begin{thebibliography}{99}
\bibitem{b1}
Celeghini E., Giachetti R., Sorace E., Tarlini M.,
J. Math. Phys {\bf 31}, 2548 (1990)
\bibitem{b2}
Celeghini E., Giachetti R., Sorace E., Tarlini M.,
J. Math. Phys {\bf 32}, 1155 (1991)
\bibitem{b3}
Celeghini E., Giachetti R., Sorace E., Tarlini M.,
J. Math. Phys {\bf 32}, 1159 (1991)
\bibitem{b4}
Celeghini E., Giachetti R., Sorace E., Tarlini M.,
{\sl ``Contractions of quantum groups" in ``Quantum groups",}\/ Lecture Notes
in Mathematics 1510, 221, Springer Verlag, (1992)
\bibitem{b5}
Lukierski J., Ruegg H., Nowicki A., Tolstoy V., Phys. Lett. {\bf B264}, 331 (1991)
\bibitem{b6}
Lukierski J., Nowicki A., Ruegg H., Phys. Lett. {\bf B293}, 344 (1993)
\bibitem{b7}
Zaugg P., J. Math. Phys. {\bf 36}, 1547 (1995)
\bibitem{b8}
Zaugg P., J. Math. Phys. {\bf 28}, 2589 (1995)
\bibitem{b9}
Sobczyk  J., Czech. J. Phys. {\bf 46}, 265 (1996)
\bibitem{b10}
Woronowicz S.L., RIMS {\bf 29}, 117 (1987)
\bibitem{b11}
Ma\'slanka P., J. Math. Phys. {\bf 35}, 1976 (1994)
\bibitem{b12}
Ballesteros A., Celeghini E., Giachetti R., Sorace E., Tarlini M.,
J. Phys. {\bf A 26}, 7495 (1993)
\bibitem{b13}
Stachura P., Lett. Math. Phys. {\bf 25}, 175 (1992)
\bibitem{b14}
Giller S., Gonera C., Kosi\'nski P., Ma\'slanka P., Acta Phys. Pol. {\bf B28}, 1121 (1997)
\bibitem{b15}
Woronowicz S.L., Commun. Math. Phys. {\bf 122}, 125 (1989)
\end{thebibliography}
\end{document}